\documentclass[11pt]{article}
\usepackage{amssymb,amsfonts,amsmath,amsthm,cite,color}
\usepackage{epsfig}
\newtheorem{thm}{Theorem}[section]

\newtheorem{prop}{Proposition}[section]
\newtheorem{lem}{Lemma}[section]
\newtheorem{exa}{Example}[section]
\newtheorem{defi}[thm]{Definition}

\makeatletter \@addtoreset{equation}{section}

\def\pf{\noindent {\it Proof.\ }}
\def\qed{\hfill \rule{4pt}{7pt}}

\begin{document}
\begin{center}
{\huge {\bf Computing  the decomposition group of a zero-dimensional
ideal by elimination method}}

\vskip 6mm

{\small Yongbin Li \\[3mm]

School of Mathematical Sciences\\
University of Electronic Science and Technology of China\\
Chengdu 611731, P.R. China \\[3mm]

yongbinli@uestc.edu.cn }
\end{center}

\noindent {\bf Abstract.} In this note, we show that the
decomposition group $Dec(I)$ of a zero-dimensional radical ideal $I$
in ${\bf K}[x_1,\ldots,x_n]$ can be represented as the direct sum of
several symmetric groups of polynomials based upon using Gr\"{o}bner
bases. The new method makes a theoretical contribution to discuss
the decomposition group of $I$  by using Computer Algebra without
considering the complexity. As one application, we also present an
approach to yield new triangular sets in computing triangular
decomposition of polynomial sets ${\mathbb P}$ if $Dec(<{\mathbb
P}>)$ is known.

\noindent \textbf{Keywords:} Decomposition group, symmetric group,
zero-dimensional ideal, Gr\"{o}bner basis, characteristic
polynomial,triangular decomposition


\section{Introduction}
Let {\bf K} be  a perfect field  and ${\bf K}[x_1,\ldots,x_n]$ (or
${\bf {K[x]}}$ for short) the ring of polynomials in
$x_1,\ldots,x_n$ with coefficients in {\bf K}.  Let $S_n$ be the
symmetric group on $n-$ elements. Consider the canonical action of
$S_n$ on a special ideal $I$ of ${\bf {K[x]}}$. We focus upon
computation of the subgroup
 $Dec(I)\leq S_n$ which leaves $I$ globally
invariant. The group $Dec(I)$ is  called the {\it decomposition
group} of $I$.

 $Dec(I)$ is up to the Galois group of $F\in  {\bf K}[x]$ where $I$ is the defining ideal of
 an irreducible polynomial $F\in {\bf K}[x]$ of degree $n$. Anai, Noro and Yokohama in \cite{HA}  give an algorithm for computing $Dec(<{\mathbb
  T}>)$  which can be executed with $O(n^4).$ Ines, Abdeljaouad
  {\it et
al.} in \cite{04:dec} (unpublished) present the other one which is
bound by $O(n^3)$ by constructing an increasing chain of subgroups.

In this paper,  we consider $I$  as a zero-dimensional radical
ideal. Some properties of $Dec(I)$ are described in Section 2. Lemma
3.2 in Section 3 shows an approach to compute an useful partition of
the set $\{1,2,\dots,n\}$ by the fine structure of ${\bf K}[{\bf
x}]/I$ using Gr\"{o}bner bases. Based upon the  partition, an
alternative method is given by Theorem 3.1 in Section 3 for
computing $Dec(I)$ without considering the complexity. Our strategy
of the computation of $Dec(I)$ is to compute the direct product of
the symmetric groups $Symm(F_k)$ of some multivariate polynomials
$F_k$.

As an application of $Dec(I)$ where $I$ is a zero-dimensional ideal,
Theorem 5.1 in Section 5 is helpful to get some new triangular sets
 in triangular decomposition of $I$ when $Dec(I)$ is given.

\section{Preliminaries}

\vspace{.1in}

 The extension field $\bf \tilde {K}$ of
{\bf K} is algebraically closed in this paper. A {\it polynomial
set} is a finite set $ {\mathbb  P}$ of nonzero polynomials in $\bf
{K[x]}$. The ideal of ${\bf {K[x]}}$ generated by all elements of $
{\mathbb  P}$ is denoted by $<{\mathbb P}>$. The  zero set of
${\mathbb P}$ in $\tilde{{\bf
 K}}^n$,   denoted by ${\rm
Zero}_{\tilde{\bf K}^{n}}({\mathbb
 P})$, is called the {\it affine variety} defined by ${\mathbb P}$.
 It is obvious that ${\rm Zero}_{\tilde{\bf K}^{n}}({\mathbb P})={\rm Zero}_{\tilde{\bf K}^{n}}(<{\mathbb
P}>).$ If ${\rm Zero}_{\tilde{\bf K}^{n}}({\mathbb
 P})$ is a finite set, then $<{\mathbb P}>$ is called a
{\it zero-dimensional} ideal of ${\bf {K[x]}}$. An ideal $I$ in
${\bf K[x]}$ such that $I =\sqrt{I} $ is called a {\it radical
ideal}.

 The next definition of
$Dec(I)$ generalizes  the one mentioned in \cite{HA,04:dec}.

\begin{defi}

 Let $I$ be a zero-dimensional radical ideal
 in ${\bf K}[x_1,\ldots,x_n]$. The following permutation group
$Dec(I)$ is called the decomposition group of $I$
$$Dec(I)\triangleq \{\sigma \in S_n|\,f(x_{\sigma(1)},\ldots,x_{\sigma(n)})\in I,\;\forall\,f\in I\}.$$
\end{defi}

\begin{exa}

 Let $I=<x_1x_2x_3+1,x_1+x_2+x_3,x_1x_2+x_2x_3+x_3x_1>,$ then
$Dec(I)=S_3.$

\end{exa}

The next subgroup of $S_n$  which is different from $Dec(I)$
concerns only one polynomial.

\begin{defi}
Let $F\in{\bf R}[t_1,\ldots,t_n]$ where ${\bf R}$ is a commutative
ring. The following permutation group $Sym(F)$ is called the
symmetric group of $F$, $$ Sym(F)\triangleq\{\sigma \in
S_n|\,F(t_{\sigma(1)},\ldots,t_{\sigma(n)})=F\}.$$

\end{defi}
\begin{exa}

 Let $F=(t_1+t_3)(t_2+t_4)\in {\bf K}[t_1,t_2,t_3,t_4],$ then
$Sym(F)=D_4.$

\end{exa}

To understand the relation between $Dec(I)$ and $Sym(F)$, we have
the following obvious claim.

\begin{lem}
If $f\in {\bf K}[{\bf x}]$, then $Sym(f)=Dec(<f>).$
\end{lem}

Referring to \cite{Cox06}, the quotient of ${\bf {K[x]}}$ modulo
$I$, written ${\bf {K[x]}}/I$, is the set of equivalence classes for
congruence modulo $I$:
$${\bf {K[x]}}/I\triangleq\{[f]\,|\,f\in {\bf {K[x]}}\}$$
which is a finitely generated ${\bf K}$-algebra.

For any $\sigma \in S_n$, we denote  $\psi_{\sigma}\,:\,{\bf
{K[x]}}\rightarrow {\bf {K[x]}}$ for the following ${\bf K}$-algebra
homomorphism given by
$$\psi_{\sigma}(g(x_1,\ldots,x_n))=g(x_{\sigma(1)},\ldots,x_{\sigma(n)})$$
for any $g\in {\bf {K[x]}}.$ It is  obvious  that $\psi_{\sigma}$ is
a ${\bf K}$-algebra isomorphism. One can easily see that the next
claim.

\begin{lem} \label{I-stable}
  With the above situation, if $\sigma \in Dec(I)$, then
$\psi_{\sigma}(I)=I.$
  \end{lem}

Based upon the above result, we proceed to discuss some interesting
properties of $Dec(I)$.

\begin{prop}
\label{DecByZero}
 Let $I$ be a zero-dimensional radical ideal in ${\bf K}[{\bf
x}]$. Then
$$Dec(I)=\{\sigma \in S_n|\,(a_{\sigma(1)},\ldots,a_{\sigma(n)})\in {\rm Zero}_{{\tilde{\bf K}}^n}(I),\;\forall\,(a_1,\ldots,a_n)\in {\rm Zero}_{{\tilde{\bf K}}^n}(I)\}.$$

\end{prop}

\pf For any $\sigma \in Dec(I)$ and $(a_1,\ldots,a_n)\in {\rm
Zero}_{{\tilde{\bf K}}^n}(I)$,
$(a_{\sigma(1)},\ldots,a_{\sigma(n)})\in {\rm Zero}_{{\tilde{\bf
K}}^n}(I)$ follows from $\psi_{\sigma}(I)=I$ by Lemma
\ref{I-stable}. Going in the other direction, considering that $I$
is a zero-dimensional radical ideal, we have $I({\rm
Zero}_{\tilde{{\bf K}^n}}(I))=I$ by the Nullstellensatz. Given a
$\sigma \in S_n$, $(a_{\sigma(1)},\ldots,a_{\sigma(n)})\in {\rm
Zero}_{{\tilde{\bf K}}^n}(I)$ for any $(a_1,\ldots,a_n)\in {\rm
Zero}_{{\tilde{\bf K}}^n}(I)$ implies that $\psi_{\sigma}(I)=I.$
Thus, our claim keeps true. \qed

\begin{prop} Let $I$ be a zero-dimensional radical ideal in ${\bf K}[{\bf x}]$. Then
$$Dec(I)=\{\sigma \in S_n|\,[f(x_{\sigma(1)},\ldots,x_{\sigma(n)})]\in {\bf {K[x]}}/I,\;\forall\,[f]\in{\bf {K[x]}}/I\}.$$

\end{prop}


\section{An alternative  method for computing $Dec(I)$}

\subsection{A direct mehtod}
Here we firstly give a direct method for computing $Dec(I)$ if ${\rm
Zero}_{{\tilde{\bf K}}^n}(I)$ is known when $I$ is a
zero-dimensional radical ideal in  ${\bf K}[{\bf x}]$.

Let  $ {\rm Zero}_{{\tilde{\bf K}}^n}(I)=\{P_1,\ldots,P_N\}$ with
$P_i=(x_1^{(i)},\cdots,x_n^{(i)})$ for $i=1,2,\dots, N.$  Now set
$S_i=\{x_i^{(k)}\,|\, k=1,2,\dots,N\}$ for $i=1,2,\dots,n.$ We
denote ${\bf S}$ for the set of all $S_i.$

\begin{exa}
Let ${\rm Zero}_{{\tilde{\bf
K}}^n}(I)=\{(2,3,5,6),(2,5,3,6),(2,5,6,3),(2,6,5,3)\}.$ Then, ${\bf
S}=\{S_1,S_2,S_3,S_4\}$ where
$$S_1=\{2\},S_2=\{3,5,6\},S_3=\{3,5,6\},S_4=\{3,6\}.$$

\end{exa}

Define two elements $ S_i,S_j\in {\bf S}$ to be $I$-equivalent,
denoted $S_i\sim_I S_ j$, if $S_i=S_ j$. One can easily check that
$\sim_I$ is an equivalence relation, here we omit the details. The
equivalence classes of $\sim_I$ yields a partition of ${\bf S}.$ The
above partition  can naturally  induce the following partition of
$\{1,2,\dots,n\}$, denoted $P_{I}$.

\begin{lem}

If ${\rm Zero}_{{\tilde{\bf K}}^n}(I)$ is given when $I$ is a
zero-dimensional radical ideal in ${\bf K}[{\bf x}],$ then  $P_I$ is
a partition of $\{1,2,\dots,n\}.$

\end{lem}

\begin{exa} Continued from the above example, we have
 $P_{I}=\{B_1,B_2,B_3\}$ with $B_1=\{1\},B_2=\{4\},B_3=\{2,3\}.$
\end{exa}

In fact, we can obtain $P_I$ by virtue of ${\bf K}[{\bf x}]/I$
without computing ${\rm Zero}_{{\tilde{\bf K}}^n}(I)$ in the next
section.

\subsection{Computing $Dec(I)$ }
When ${\rm Zero}_{{\tilde{\bf K}}^n}(I)$  is finite,  ${\bf K}[{\bf
x}]/I$ is a finite-dimensional vector space over ${\bf K}$,which has
a nice structure. The next results due to Chapter 2 of \cite{Cox06}.
\begin{prop}
The map $\Phi  : K[{\bf x}]/I\rightarrow End_{K}(K[{\bf x}]/I)$
defined by $\Phi([f])=m_f$ where $m_f$ is a $K$-linear mapping from
$K[{\bf x}]/I$ to $K[{\bf x}]/I$ defined by $m_f([g])=[fg]$ is a
$K$-algebra homomorphism.
\end{prop}

\begin{prop}
\label{f(P)}
 With the above situation, let ${\rm Zero}_{{\tilde{\bf K}}^n}(I)
=\{P_1,\ldots,P_N\}$.
 For any $f\in K[{\bf x}]$, considering $m_f\in End_{\bar{K}}(\bar{K}[{\bf
x}]/I)$, then $f(P_i)$ is an eigenvalue of $m_f$ for $i=1,\dots,N.$

\end{prop}

In our method, the partition  $P_{I}$ of $\{1,2,\dots,n\}$ plays a
crucial role. Next we will obtain $P_{I}$ only by $m_{x_i}.$

\begin{lem} \label{partition}
With the above situation, $S_i\sim_I S_ j$ if and only if $f_i=f_j$
where $f_k$ is the character polynomial of $m_{x_i}$ for $1\leq
k\leq n.$
\end{lem}
\pf Let  $ {\rm Zero}_{{\tilde{\bf K}}^n}(I)=\{P_1,\ldots,P_N\}$
with $P_i=(x_1^{(i)},\cdots,x_n^{(i)})$ for $1\leq i \leq N.$ With
the same notation, $S_i=\{x_i^{(k)}\,|\, k=1,2,\dots,N\}$ for
$i=1,2,\dots,n.$ By Proposition \ref{f(P)}, we know that the roots
of $f_i$ over $\tilde{{\bf K}}$ are all elements of  $S_i.$
Therefore our claim holds. \qed

\begin{exa}
Let $I=<x_1x_2,x_2+x_1+1,x_3(x_3+1),x_3+x_4-1>.$ One will compute
that
$$ f_1=f_2=f_3=\lambda^4+2\lambda^3+\lambda^2,\;
f_4=\lambda^4-6\lambda^3+13\lambda^2-12\lambda+4
$$
where  $f_i(\lambda)$ is the character polynomial of $\,m_{x_i}\in
End_{\bar{K}}(\bar{K}[{\bf x}]/I)$ for $1\leq i \leq 4$, see later
in the section. Thus,$$P_{I}=\{\{1,2,3\},\{4\}\}.$$

\end{exa}

\begin{thm} \label{main}
With the above situation, let the partition $P_{I}=\{ B_k\subseteq
\{1,2,\dots,n\}\;|\;k=1,\dots,s \}$ with $B_k=\{n_1,\ldots,n_k\}.$
Then the decomposition group of $I$ is the direct product of the
following subgroups $Sym(F_k)$,
 $$ Dec(I)=\prod_{k=1}^s Sym(F_k)$$
where $F_k\in {\bf K}(\lambda)[t_1,\ldots,t_n] $ is  the character
polynomial $t_{n_1}m_{x_{n_1}}+\cdots+t_{n_k}m_{x_{n_k}}$ with new
variables $t_1,\ldots,t_n.$

\end{thm}
\pf Consider $P_{I}=\{ B_k\subseteq \{1,2,\dots,n\}\;|\;k=1,\dots,s
\}$ with $B_k=\{n_1,\ldots,n_k\},$

For any $\sigma \in Dec(I)$, it follows from Lemma \ref{I-stable}
that $\sigma (i)\in B_k$ for any $i\in B_k$ with $k=1,\dots,s.$ This
implies that
$$t_{n_1}m_{x_{n_1}}+\cdots+t_{n_k}m_{x_{n_k}}=t_{\sigma(n_1)}m_{x_{\sigma(n_1)}}+\cdots+t_{\sigma(n_k)}m_{x_{\sigma(n_k)}}$$
for any $1\leq k \leq s.$ Therefore,
$F_k(t_{\sigma(n_1)},\ldots,t_{\sigma(n_k)})=F_k(t_{n_1},\ldots,t_{n_k})$
by Proposition \ref{f(P)}. So, $$\sigma \in \prod_{k=1}^s
Sym(F_k).$$

\qed

\begin{exa}
Continued from the above example,
 $Dec(I)=Sym(F_1)\times Sym(F_2)=\{(1),(1,2)\}\times \{(4)\}=\{(1),(1,2)\}$  where
\begin{eqnarray*}
 F_1&=&\left( \lambda+t_{{2}} \right)  \left( \lambda+t_{{2}}+t_{{3}} \right)  \left( \lambda+t_{{1}} \right)  \left( \lambda+t_{{3}}+t_{{1}}
 \right),\\
  F_2&=&\left( \lambda-t_{{4}} \right) ^{2} \left( \lambda-2\,t_{{4}} \right) ^{2}.
\end{eqnarray*}

\end{exa}

{\bf An alternative method :} Given a zero-dimensional radical ideal
$I$ in ${\bf K}[{\bf x}]$, we try to compute $Dec(I)$ as follows.
\begin{description}
 \item [Step 1.] Compute a Gr\"obner basis of $I$.
\item[Step 2.]  Compute each $m_{x_i}$ for $1\leq i\leq n.$
  \item [Step 3.] Compute each character polynomial $f_i$ of  $m_{x_i}$ for $1\leq i\leq n.$
 \item [Step 4.] Construct the partition $P_{I}$ of $\{1,2,\dots,n\}$
 by Lemma \ref{partition}.
 \item [Step 5.] Compute the $Dec(I)$ by Theorem \ref{main}.
  \end{description}

\begin{exa}

 Consider
$I=<x_1x_2,x_2+x_1+1,x_3(x_3+1),x_3+x_4-1>.$
 \begin{description}

\item[\textbf{Step 1.}] Compute a Gr\"obner basis $G$
of $I$ as follows,

$G=[{x_{{4}}}^{2}-3\,x_{{4}}+2,x_{{4}}+x_{{3}}-1,{x_{{2}}}^{2}+x_{{2}},x_
{{1}}+x_{{2}}+1].$

\item[\textbf{Step 2.}] Compute each $m_{x_i}$ for $1\leq i\leq 4:$
$$m_{x_1}=\left[ \begin {array}{cccc}
-1&0&-1&0\\\noalign{\medskip}0&-1&0&-1\\\noalign{\medskip}0&0&0&0\\\noalign{\medskip}0&0&0&0\end
{array}
 \right], m_{x_2}=\left[ \begin {array}{cccc} 0&0&1&0\\\noalign{\medskip}0&0&0&1
\\\noalign{\medskip}0&0&-1&0\\\noalign{\medskip}0&0&0&-1\end {array}
 \right], $$

$$m_{x_3}=\left[ \begin {array}{cccc}
1&-1&0&0\\\noalign{\medskip}2&-2&0&0\\\noalign{\medskip}0&0&1&-1\\\noalign{\medskip}0&0&2&-2\end
{array}
 \right], m_{x_4}= \left[ \begin {array}{cccc} 0&1&0&0\\\noalign{\medskip}-2&3&0&0
\\\noalign{\medskip}0&0&0&1\\\noalign{\medskip}0&0&-2&3\end {array}
 \right].$$

\item[\textbf{Step 3.}] Compute each character polynomial $f_i$ of  $m_{x_i}:$

$$f_1=f_2=f_3=\lambda^4+2\lambda^3+\lambda^2,$$
$$f_4=\lambda^4-6\lambda^3+13\lambda^2-12\lambda+4.$$

\item[\textbf{Step 4.}] Construct the partition $P_{I}=\{\{1,2,3\},\{4\}\}$ of $\{1,2,\dots,4\}$
 by Lemma \ref{partition}.

\end{description}

\end{exa}


\begin{description}

\item[\textbf{Step 5.}] Compute each character
polynomial $F_k $,
\begin{eqnarray*}
F_1&=& |t_1m_{x_1}+t_2m_{x_2}+t_3m_{x_3}-\lambda I|\\
&=&\left( \lambda+t_{{2}} \right) \left( \lambda+t_{{2}}+t_{{3}}
\right) \left( \lambda+t_{{1}} \right) \left(
\lambda+t_{{3}}+t_{{1}}
 \right),\\
  F_2&=&|t_4m_{x_4}-\lambda I|=\left( \lambda-t_{{4}} \right) ^{2} \left( \lambda-2\,t_{{4}} \right) ^{2}.
  \end{eqnarray*}

 Compute $Sym(F_1)= \{(1),(1,2)\}$ and  $Sym(F_2)= \{(4)\}.$
 \vskip 5pt
Thus, $Dec(I)=Sym(F_1)\times Sym(F_2)=\{(1),(1,2)\}$ by Theorem
\ref{main}.

\end{description}

\vspace{.1in}

\section{One application}

Throughout this section, with the same notation in \cite{ChenM}, we
denote ${\bf K}[x_1<\ldots <x_n]$ by the ring ${\bf K}[{\bf x}]$
with ordered variables $x_1<\ldots <x_n.$ A polynomial ordered set $
{\mathbb T}=[f_1,\ldots,f_n]$ which could be written in the form
  \begin{eqnarray*}
  \label{Tr}
f_1\in {\bf K}[x_1]\setminus {\bf K},\;f_i=f_i(x_1,\ldots,x_{i})\in
{\bf K}[x_1,x_2,\ldots,x_i]\setminus {\bf K}[x_1,x_2,\dots,x_{i-1}]
\end{eqnarray*}
for $2\leq i \leq n$ is called a  {\it zero-dimensional triangular
set} in ${\bf K}[x_1<\ldots <x_n]$. For $f_j\in {\mathbb T}$,  ${\rm
lcoeff}(f_j,x_j)$ is called the {\it leading coefficient} of $f_j$,
denoted by ${\rm lc}(f_j)$.
 A zero-dimensional triangular
set ${\mathbb T}=[f_1,\ldots,f_n]$ is call a {\it regular set} or
{\it chain} if ${\rm res}({\rm lc}(f_j),{\mathbb T})\ne 0$  where
${\rm res}({\rm lc}(f_j),{\mathbb T})$ stands for the {\it
resultant} of ${\rm lc}(f_j)$ with respect to ${\mathbb { T}}$ for
$j=2,\cdots,s.$ See \cite{aub,Kalk:regular,gc:regular,yz:RegularSet,
wd:regular} for the details.

Given a polynomial set $ {\mathbb  P}$ such that $I=<{\mathbb P}>$
is a zero-dimensional radical ideal in $\bf {K[x]},$ based upon one
of several successful algorithms for triangular decomposition
presented in
\cite{Kalk:regular,lazard91,gd:book,wd:book,yz:book,wu78,wu:book} ,
one may get the following zero decomposition of the form
$${\rm Zero}_{\tilde{\bf K}^{n}}(I)= \bigcup^e_{i=1}{\rm Zero}_{\tilde{\bf K}^{n}}({\mathbb T}_i).$$
For convenience, we assume that  each ${\mathbb T}_i$ in the above
decomposition is the regular set or chain  in ${\bf K}[x_1<\ldots
<x_n]$.

As one application of the decomposition group $Dec(I)$, the
following result is helpful to compute new triangular sets based on
a known triangular set in the above algorithms for triangular
decomposition when $Dec(I)$ is given.

\begin{thm}
Let ${\bf K}$ be a field, let $I$ be a zero-dimensional  ideal, and
let ${\mathbb T}$ be  a regular set in ${\bf K}[x_1<\cdots<x_n]$. If
$\psi_{\sigma}({\mathbb T})$ is  a triangular set, then
  $${\rm Zero}_{\tilde{\bf K}^{n}}(\psi_{\sigma}({\mathbb T}))\subseteq {\rm Zero}_{\tilde{\bf K}^{n}}(I).$$

\end{thm}

\pf Consider ${\rm Zero}_{\tilde{\bf K}^{n}}({\mathbb T})\subseteq
{\rm Zero}_{\tilde{\bf K}^{n}}(I)$, it follows from Proposition
\ref{DecByZero} that $${\rm Zero}_{\tilde{\bf
K}^{n}}(\psi_{\sigma}({\mathbb T}))\subseteq {\rm Zero}_{\tilde{\bf
K}^{n}}(I).$$ \qed

The next example wants to illustrate how triangular decomposition
becomes simpler at some time by applying the above result.

\begin{exa}
Let $I=<f_1,\ldots,f_5>$ where
\begin{eqnarray*}
f_1&=&x_{{1}}+x_{{2}}+x_{{3}}+x_{{4}}+x_{{5}},\\
f_2&=&x_{{1}}x_{{2}}+x_{{2}}x_{{3}}+x_{{3}}x_{{4}}+x_{{4}}x_{{5}}+x_{{5}}x_{
{1}},\\
f_3&=&x_{{1}}x_{{2}}x_{{3}}+x_{{3}}x_{{2}}x_{{4}}+x_{{3}}x_{{5}}x_{{4}}+x_{{
5}}x_{{1}}x_{{4}}+x_{{1}}x_{{2}}x_{{5}},\\
f_4&=&x_{{3}}x_{{2}}x_{{1}}x_{{4}}+x_{{4}}x_{{3}}x_{{5}}x_{{2}}+x_{{3}}x_{{4
}}x_{{5}}x_{{1}}+x_{{5}}x_{{2}}x_{{1}}x_{{4}}+x_{{1}}x_{{2}}x_{{5}}x_{
{3}},\\
f_5&=&x_{{5}}x_{{3}}x_{{2}}x_{{1}}x_{{4}}+1.
\end{eqnarray*}
It is obvious that $Dec(\sqrt{I})=S_5.$

Using the Triangularize command of the library of RegularChains in
Maple, see \cite{regularchain} for the details, we can compute 14
regular sets ${\mathbb T}_1,\ldots,{\mathbb T}_{14}$ such that
$${\rm Zero}_{\tilde{\bf K}^{n}}(\{f_1,\ldots,f_5\})= \bigcup^{14}_{i=1}{\rm Zero}_{\tilde{\bf K}^{n}}({\mathbb T}_i),$$
where
\begin{eqnarray*}
{\mathbb
T}_1&=&[x_{{1}}+1,x_{{2}}+1,1+x_{{3}},1-3\,x_{{4}}+{x^2_{{4}}},-3+x_{{4}}+x
_{{5}}],\\
 {\mathbb
T}_2&=&[{x^2_{{1}}}-3\,x_{{1}}+1,x_{{2}}+1,1+x_{{3}},1+x_{{4}},x_{{1}}-3+x_
{{5}}]
,\\
{\mathbb
T}_3&=&[x_{{1}}+1,1-3\,x_{{2}}+{x^2_{{2}}},1+x_{{3}},1+x_{{4}},-3+x_{{2}}+x
_{{5}}],\\
{\mathbb
T}_4&=&[x_{{1}}+1,x_{{2}}+1,-3\,x_{{3}}+1+x^2_{{3}},1+x_{{4}},-3+x_{{5}}+
x_{{3}}].
 \end{eqnarray*}

The following fact shows that all   ${\mathbb T}_2,{\mathbb T}_2$
and ${\mathbb T}_4$ can be obtained by ${\mathbb T}_1.$

\begin{eqnarray*}
{\mathbb T}_2=\psi_{(1,4)}({\mathbb
T}_1)&=&[{x^2_{{1}}}-3\,x_{{1}}+1,x_{{2}}+1,1+x_{{3}},1+x_{{4}},x_{{1}}-3+x_
{{5}}]
,\\
{\mathbb T}_3=\psi_{(2,4)}({\mathbb
T}_1)&=&[x_{{1}}+1,1-3\,x_{{2}}+{x^2_{{2}}},1+x_{{3}},1+x_{{4}},-3+x_{{2}}+x
_{{5}}],\\
{\mathbb T}_4=\psi_{(3,4)}({\mathbb
T}_1)&=&[x_{{1}}+1,x_{{2}}+1,-3\,x_{{3}}+1+x^2_{{3}},1+x_{{4}},-3+x_{{5}}+
x_{{3}}].
 \end{eqnarray*}

\end{exa}

\end{document}